\documentclass[12pt,a4paper]{article}
\usepackage[T2A]{fontenc}
\usepackage[english,russian]{babel}

\fontsize{12}{21} \selectfont

\newtheorem{theorem}{Theorem}[section]
\newtheorem{lemma}{Lemma}[section]
\newtheorem{note}{Note}[section]

\newtheorem{definition}{Definition}[section]

\begin{document}

\noindent
\bigskip
\begin{center}
{\large\bf
The analog of the Schauder inequality for closed surfaces in Euclidean spaces\\
}
\end{center}
\medskip

\begin{center}
{\bf Andrei~I.~Bodrenko}
\footnote{ Andrei~I.~Bodrenko, associate professor, Department of Mathematics,
Volgograd State University,\\
Universitó Prospekt 100, Volgograd, 400062, RUSSIA.
Email: bodrenko@mail.ru}
\end{center}

\begin{center}
{\bf Abstract}\\
\end{center}

{\small
The analog of the Schauder inequality for closed surfaces in Euclidean spaces
is obtained in this article.
}

\section*{Introduction}

Let $E^{n}$ be $n$-dimensional ($n>1$) Euclidean space.
Let $D$ be the finite domain in $E^n$, $\partial D$ be the boundary of $D,$
$\bar D$ be the closure of $D.$
Let $(x^{1},...,x^{n})$ be the Cartesian coordinates in $E^{n}.$

\begin{definition}
\label{definition 1}.
We say, that function $f$ on $D$ is of class
$C^{k,s}(\bar D)$ , $k\geq 0,$
$s\in (0,1),$ if it has continuous partial derivatives
up to $k$-th order inclusively and bounded value
$$
|f|_{(D) k,s}= \sum_{|i|=0}^{k} \sup_{x\in D} | \partial^{|i|} f(x)|
+
\sum_{|i|=k} \sup_{x_1, x_2 \in D} \frac{|\partial^{|i|} f(x_{1})-
\partial^{|i|}f(x_{2})|}
{|x_{1}-x_{2}|^{s}}.
\eqno(1)
$$
 Partial derivatives of function $f(x)$ are denoted by
 $\partial^{|i|} f(x)\equiv$
$\frac{\partial^{|i|}f(x)}
{\partial^{i_{1}}x^{1}...\partial^{i_{n}}x^{n}}$,
where $|i|=i_{1}+ ... +i_{n}$ is the order of derivative.
$|x|=(\sum_{i=1}^{n}(x^{i})^{2})^{1/2}$, where $(x^{1},...,x^{n})$ are
the coordinates of point $x\in E^{n}$.
\end{definition}

The value $|f|_{(D) k,s}$ we call
the norm of function $f$  in the space $C^{k,s}(\bar D)$.
It is known (see [1]) that the space $C^{k,s}(\bar D)$
with norm denoted by formula $(1)$ is complete normed space.

We define the cylinder $C_{R,L}$ in $E^{n}$ by the following formula:
$$C_{R,L}=\Biggl\{x: \sum_{i=1}^{n-1}(x^{i})^{2}
< R^{2}, -2LR < x^{n} < 2LR \Biggr\},$$
where $ L=const>0, R=const>0, $
and let $x=(0,0,...,0)$ be called its center.

\begin{definition}
\label{definition 2}.
 Domain $D$ is called strictly Lipschitzian if
for every point $x_{0}\in \partial D$ it can be introduced the coordinates
$$
u^{k}=\sum_{l=1}^{n} c^{k}_{l}(x^{l}-x^{l}_{0}),
\qquad
k=\overline {1,n},
$$
where $||c^{k}_{l}||$ is orthogonal
matrix such that the intersection of $\partial D$ and the cylinder
$\bar C_{R,L}$
corresponding to the coordinates $\{u^{k}\}$ , is given by the equation
$$
u^{n}=\omega (u^{,n}),
\qquad
u^{,n}\equiv (u^{1},...,u^{n-1}),
$$
where $\omega (u^{,n})$ is Lipschitzian function for $|u^{,n}|\leq R$
with Lipschitz constant bounded by $L,$
and
$$
\bar D \cap \bar C_{R,L} =\Biggl\{ u: |u^{,n}|\leq R,
\quad
-2LR\leq u^{n}\leq\omega (u^{,n}) \Biggr\}.
$$
\end{definition}
The numbers $R$ and $L$ are fixed for the domain $D$.

Arbitrary convex domain is strictly Lipschitzian (see [1]).

Let $x_{0}=(x^{1}_{0},...,x^{n}_{0})$ be a point of $\partial D$,
where surface $\partial D$ has tangent plane.

\begin{definition}
\label{definition 3}.
We call $(u^{1},...,u^{n})$ the specific local coordinate system
with origin at the point $x_{0}$ if the coordinates $\{u^{k}\}$ and
$\{x^{k}\}$
satisfy the following equations: $u^{k}=\sum_{l=1}^{n} c^{k}_{l}(x^{l}-x^{l}_{0}),$
$ k=\overline {1,n},$
and the axis $u^{n}$ is directed to the normal of $\partial D$
at the point $x_{0}$, that is outward for $D$.
\end{definition}

\begin{definition}
\label{definition 4}.
Domain $D$ is called the domain of class $C^{l,s}$, $l\geq 1$,
if it is strictly Lipschitzian and the coordinates $\{u^{k}\}$
that are given in the definition 2
are the specific local coordinates, the function $u^{n}=\omega (u^{,n})$,
defining the equation of the surface $\partial D$
is of class $C^{l,s}(|u^{,n}|\leq R)$.
\end{definition}

\begin{definition}
\label{definition 5}.
 We will say that the boundary $\partial D$ of domain $D$
is of class $C^{l,s}$
if for every point $x_{0}\in \partial D$ there can be introduced
the the specific local coordinates such that
the function $u^{n}=\omega (u^{,n})$ is of class
$C^{l,s}(|u^{,n}|\leq R).$
\end{definition}

Let the boundary $\partial D$ of the domain $D$ is of
class $C^{l_{1},s_{1}},$ $s_{1}\in (0;1)$.
Let on $\partial D$ be given the function $\varphi (x), x\in \partial D$.

\begin{definition}
\label{definition 6}.
  We will say that function $\varphi (x)$ is of
class $C^{l,s}(\partial D)$ if it as a function
of the specific local coordinates $u^{,n}=(u^{1},...,u^{n-1})$
introduced for every
point $x_{0}\in\partial D$ is of class $C^{l,s}(|u^{,n}|\leq R),$
where $|u^{,n}|\leq R$ is the base of cylinder corresponding to
the point $x_{0}.$
\end{definition}

\begin{definition}
\label{definition 7}.
  Norm $|\varphi|_{(\partial D) l,s}$ of function $\varphi$ ,
given on the surface $\partial D$ is called the greatest of the norms
$|\varphi(u^{,n})|_{(|u^{,n}|\leq R) l,s} ,$  calculated for all
points $x_{0}\in \partial D.$
\end{definition}

Let $F$ be the closed orientable hypersurface in Euclidean space $E^{n+1}.$
Let $(y^{1},...,y^{n+1})$ be the Cartesian coordinates in $E^{n+1}.$
Let $U$ be arbitrary open set on $F.$

\begin{definition}
\label{definition 8}.
  Couple $(U,h)$ is called the admissible map of class
$C^{k,s}$ if:

1) $h$ is the homeomorphism $U$ on the open ball $K_{r}$ of radius $r>0$
in $E^{n};$

2) the inverse mapping $\bar h^{-1}(x)
\equiv (f^{1}(x),...,f^{n+1}(x))$, $x\in K_{r}$,
satisfies the condition:
$f^{\alpha}\in C^{k,s}(\bar K_{r})$, $\alpha=\overline{1,n+1}$.
\end{definition}

\begin{definition}
\label{definition 9}.
  $F$ is called the surface of class $C^{k,s}$
if the following conditions hold:

1) $F$ is the surface of class $C^{k}$;

2) on $F,$ there exists the finite
 aggregate $\{(U_{i},h_{i})\}_{i=\overline {1,N}}$ of admissible
maps of class $C^{k,s},$ where the collection of sets
$(U_{i})_{i=\overline {1,N}}$ is open covering of $F$;

3) if $U_{i}\cap U_{j}\not =\emptyset$ then
the mapping $h_{j}\circ h^{-1}_{i}$
is diffeomorphism of class $C^k$ of
$h_{i}(U_{i}\cap U_{j})$ on set $h_{j}(U_{i}\cap U_{j})$.
\end{definition}

\begin{definition}
\label{definition 10}.
The aggregate $\{(U_{i},h_{i})\}_{i=\overline {1,N}}$ of the admissible maps
considered in definition 9 is called the admissible atlas of class
$C^{k,s}$ of hypersurface $F.$
\end{definition}

\begin{definition}
\label{definition 11}.
Function $f$ determined on surface $F,$
is called the function of class $C^{p,s}(F)$, $p<k $, if

1)  on hypersurface $F,$ there exists the admissible
 atlas $\{(U_{i},h_{i})\}_{i=\overline {1,N}}$ of class $C^{k,s}$ and

2) $f\circ h_{i}^{-1}\in C^{p,s}(\bar K_{r_{i}})$, $i=\overline {1,N}$.
\end{definition}

We fix the admissible atlas $\{(U_{i},h_{i})\}_{i=\overline {1,N}}$
of class $C^{k,s}$ of hypersurface $F$.

Let be given the function $f\in C^{p,s}(F)$.
The norm of function $f$ in space $C^{p,s}(F)$ is determined
by the following formula:
$$|f|_{(F) p,s}= \max_{i} |f\circ h^{-1}_{i}|_{(K_{r_{i}}) p,s}.$$

We will prove that the obtained normed space is complete.
Let $\{f_{m}\}_{m=1}^{\infty}$ be the Cauchy sequence of functions $f_m$
of class $C^{p,s}(F),$ therefore $\forall \varepsilon>0 $ and
for every natural number $l$

$$|f_{m+l}-f_{m}|_{(F) p,s}< \varepsilon, $$ for all safficiently large $m.$
Then, from the definition of norm, we obtain
$$|f_{m+l}\circ h^{-1}_{i}-f_{m}\circ h^{-1}_{i}|_{(K_{r_{i}})p,s}<\varepsilon.$$
Since the function space $C^{p,s}(\bar K_{r_{i}})$ is
complete normed space, then the sequence of functions
$\{f_{m}\circ h^{-1}_{i}\}$ on $\bar K_{r_{i}}$ has limit:
$f\circ h^{-1}_{i}=\lim_{m\rightarrow\infty}f_{m}\circ h^{-1}_{i}.$
Since $f\circ h^{-1}_{i}\in C^{p,s}(\bar K_{r_{i}}),
\forall i=\overline{1,N},$ then $f\in C^{p,s}(F)$.

\section{Statement of the result.}

Let $F\in C^{3,s},$ where $s\in (0;1).$
Let $\{(U_{i},h_{i})\}_{i=\overline {1,N}}$ be the admissible
 atlas $F$ of class
$C^{3,s}$. Let $(U,h)$ be an arbitrary map from the collection
$\{(U_{i},h_{i})\}_{i=\overline {1,N}}$.
Then the hypersurface $F$ on map $(U,h)$
is determined by the following equation system:
$$y^{\alpha}\equiv h^{-1 \alpha}(x)=f^{\alpha}(x^{1},...,x^{n}),
\alpha=\overline{1,n+1}, (x^{1},...,x^{n})\in K_{r} \eqno(2)$$

Consider the differential operator $A$ on $F$
that, on map $(U,h),$ is defined by:
$$A=\sum_{k=1}^{n}\sum_{p=1}^{n}a^{kp}
\partial_{k}(\partial_{p}) + \sum_{j=1}^{n}b^{j}\partial_{j}+c. \eqno(3)$$
Let $a^{kp}=a^{pk},$ and the operator $A$ is strictly elliptic on $F,$
 i. e.
$$
\sum_{k=1}^{n}\sum_{p=1}^{n} a^{kp}(x) \zeta_{k}\zeta_{p} \geq
\nu \sum_{k=1}^{n}(\zeta_{k})^{2},
\quad
\nu=const>0,
\quad
\forall \zeta_{k},
\quad
k=\overline {1,n}.
$$

Let, on $F,$ be given a function $f$ of class $C^{2,s}(F).$ Then we have:
$$A(f\circ h^{-1}(x))=\sum_{k=1}^{n}\sum_{p=1}^{n}a^{kp}(x)
\partial_{k}(\partial_{p}(f\circ h^{-1}(x)))
+ \sum_{j=1}^{n}b^{j}(x)\partial_{j}(f\circ h^{-1}(x)) +
c(x) f\circ h^{-1}(x), $$
where $x\in K_{r}$ for every admissible map $(U,h)$ of class $C^{3,s}.$

\begin{theorem}
\label{theorem1}.
  Let function $f$ be a solution of class $C^{2,s}(F)$
of the problem: $Af=\gamma,$ where $c(x)\neq 0$ on $F,$
 $\gamma\in C^{0,s}(F),$
 $a^{kp}\in C^{0,s}(F),$ $b^{j}\in C^{0,s}(F),$ $c \in C^{0,s}(F).$
Then the following inequality holds:
$$|f|_{(F) 2,s}\leq M|\gamma|_{(F)0,s},$$ where the constant $M$ depends
on $s,n,$ the surface $F,$ the coefficients $a^{kp}, b^{j} , c $
$ ( k,p,j=\overline {1,n} )$
 and the admissible atlas on $F$
$\{(U_{i},h_{i})\}_{i=\overline {1,N}}$ of class $C^{3,s}$.
\end{theorem}

\section{Auxiliary conjectures.}

\begin{note}
\label{note1}:
Let $f$ be the function of class $C^{2,s}(F),$
$\{U_{i},h_{i}\}$ be the admissible map on $F,$
 $h_{i}(U_{i})=K_{r_{i}}$.
 Let $x_{0}\in \partial K_{r_{i}}$.
 Consider the specific local coordinates $(x^{1},...,x^{n})$
 for the point $x_{0}$ where $ x^{n}=\omega(x^{1},...,x^{n-1}).$
Then the intersection of the cylinder
$\bar{C}_{R,L}$ at the point $x_{0}$ and the surface $\partial K_{r_{i}}$
is given by:
$ x^{n}=\omega(x^{1},...,x^{n-1}), |(x^{1},...,x^{n-1})|\leq R.$
We assume that the specific local coordinates were introduced in the ball $K_{r_{i}}$.
\end{note}

Let $O_{i}(x_{0})$ be an open domain in $K_{r_{i}}$ such that
$\bar{O}_{i}(x_{0})\supset (\partial K_{r_{i}}\cap \bar{C}_{R,L})$.
Let $B_{\rho_{j}}(x_{0})$ be an open ball â $E^{n}$ of radius $\rho_{j}$
 with center at the point $x_{0}.$

We will prove several lemmas before the theorem 1.

\begin{lemma}
\label{lemma1}.
  There exist numbers $R,L,Q$ and set collection
  $\{O_{i}(x_{0})\}_{i=\overline {1,N}}$
such that for every point $x_{0}\in \partial K_{r_{i}}$ and for all
$i=\overline {1,N}$ the following conditions hold:

1) $\exists j\not =i$ such that $h_{i}^{-1}(O_{i}(x_{0}))\subset U_{j}$.

2) $\exists \rho_{j}>0$ such that
$h_{j}(h_{i}^{-1}(O_{i}(x_{0})))\subset B_{\rho_{j}}(x_{0})
\subset K_{r_{j}}$,
where $dist(B_{\rho_{j}}(x_{0}),\partial K_{r_{j}})\geq Q > 0$.
\end{lemma}

Proof of lemma 1 follows from definition 10,
compactness of $\partial K_{r_{j}},$ definition 2 and finiteness of
covering $\{U_{i}\}_{i=\overline {1,N}},$
for sufficiently small numbers $R$ and $L.$

We fix numbers $R,L,Q,$ point $x_{0}\in \partial K_{r_{i}}$
and set collection $\{O_{i}(x_{0})\}_{i=\overline {1,N}},$
that satisfy lemma 1.

\begin{lemma}
\label{lemma2}.
 Under the conditions of lemma 1, the following inequality holds:
$$
\sup_{|(x^{1},...,x^{n-1})|\leq R}
\left|
f\circ h^{-1}_{i}
(x^{1},...,x^{n-1},\omega(x^{1},...,x^{n-1}))
\right|
\leq
\sup_{u\in B_{\rho_{j}}}
\left|
f\circ h^{-1}_{j}(u)
\right|.
$$
\end{lemma}

{\bf Proof.} Since there exists a diffeomorphism of class $C^{3}$ of
the neighborhood $O_{i}(x_{0})$ into the ball $B_{\rho_{j}}(x_{0}),$
then there exist mappings:
$u^{p}=k^{p}(x^{1},...,x^{n}),p=\overline {1,n},$
$\forall x\in
O_{i}(x_{0}),$ $x^{p}=g^{p}(u^{1},...,u^{n}),  \forall u\in
h_{j}(h^{-1}_{i}(O_{i}(x_{0})).$
Hence for every point $x\in O_{i}(x_{0})$
there exists point $u\in B_{\rho_{j}}(x_{0})$
such that $f\circ h^{-1}_{i}(x)=f\circ h^{-1}_{j}(u)$.
Lemma is proved.

\begin{lemma}
\label{lemma3}.
The following inequality holds:
$$
\sup_{|(x^{1},...,x^{n-1})|\leq R}
\left|
\frac{\partial}{\partial x^{k}}
(f\circ h^{-1}_{i}
(x^{1},...,x^{n-1},\omega(x^{1},...,x^{n-1})))
\right|
\leq
$$
$$
\leq M \max_{p=\overline{1,n}}\sup_{u\in B_{\rho_{j}}}
\left|
\frac{\partial}{\partial u^{p}}(f\circ h^{-1}_{j}(u))
\right|
,
\quad
k=\overline{1,n-1},
$$
where $M=\mbox{const}>0$.
\end{lemma}

{\bf Proof.}
We have $$
\frac{\partial}{\partial x^{k}}
\left(
f\circ h^{-1}_{i}
(x^{1},...,x^{n-1},\omega(x^{1},...,x^{n-1}))
\right)
=
$$
$$
=
\frac{\partial}{\partial x^{k}}
\left(
f\circ h^{-1}_{j}(u^{1},...,u^{n})
\right)
=
\frac{\partial}{\partial u^{p}}
\left(
f\circ h^{-1}_{j}(u^{1},...,u^{n})
\right)
\frac{\partial k^{p}}{\partial x^{k}},
$$
where the point $(u^{1},...,u^{n})\in B_{\rho_{j}}(x_{0})$,
$u^{l}=k^{l}(x^{1},...,x^{n-1},\omega(x^{1},...,x^{n-1}))$,
$l=\overline{1,n}$.
Since the functions $k^{p}\in C^{3}$ then the functions
$\frac{\partial k^{p}}{\partial x^{k}}$ are bounded on
$|(x^{1},...,x^{n-1})|\leq R$.
Lemma 3 is proved.

\begin{lemma}
\label{lemma4}.
 The following inequality holds:
$$
\sup_{|(x^{1},...,x^{n-1})|\leq R}
\left|
\frac{\partial^{2}}{\partial x^{k}\partial x^{q}}
(f\circ h^{-1}_{i}
(x^{1},...,x^{n-1},\omega(x^{1},...,x^{n-1})))
\right|
\leq
$$
$$
\leq
M_{1}
\left(
\max_{p=\overline{1,n}}\sup_{u\in B_{\rho_{j}}}
\left|
\frac{\partial}{\partial u^{p}}(f\circ h^{-1}_{j}(u))
\right|
+\max_{p,r=\overline{1,n}}\sup_{u\in B_{\rho_{j}}}
\left|
\frac{\partial^{2}}{\partial u^{p}\partial u^{r}}(f\circ h^{-1}_{j}(u))
\right|
\right)
,
$$
$$
k, q= \overline{1,n-1},
$$
where $M_1=\mbox{const}>0$.
\end{lemma}

{\bf Proof.}
We have
$$
\frac{\partial^{2}}{\partial x^{k}\partial x^{q}}
\left(
f\circ h^{-1}_{i}
(x^{1},...,x^{n-1},\omega(x^{1},...,x^{n-1}))
\right)
=
$$
$$
=
\frac{\partial^{2}}{\partial u^{p}\partial u^{r}}
\left(
f\circ h^{-1}_{j}(u^{1},...,u^{n})
\right)
\frac{\partial k^{p}}{\partial x^{k}}
\frac{\partial k^{r}}{\partial x^{q}}+
\frac{\partial}{\partial u^{p}}
\left(
f\circ h^{-1}_{j}(u^{1},...,u^{n})
\right)
\frac{\partial^{2}k^{p}}{\partial x^{k}\partial x^{q}},
$$
where the point $(u^{1},...,u^{n})\in B_{\rho_{j}}(x_{0})$,
$u^{l}=k^{l}(x^{1},...,x^{n-1},\omega(x^{1},...,x^{n-1}))$,
$l=\overline{1,n}$.

Since the functions $k^{p}\in C^{3}$ therefore the functions
$\frac{\partial k^{p}}{\partial x^{k}}$ and
$\frac{\partial^{2}k^{p}}{\partial x^{k}\partial x^{q}}$
are bounded on $|(x^{1},...,x^{n-1})|\leq R$.
Therefore we obtain the proof of lemma 4.

\begin{lemma}
\label{lemma5}.
  The following inequality holds:
$$
\left|
f\circ h^{-1}_{i}(x^{1},...,x^{n-1},\omega(x^{1},...,x^{n-1}))
\right|
_{(|(x^{1},...,x^{n-1})|\leq R) 2,s}
\leq
$$
$$
\leq
M_{1}
\left(
|f\circ h^{-1}_{j}(u)|_{B_{\rho_{j}} 2,s}+
|f\circ h^{-1}_{j}(u)|_{B_{\rho_{j}} 1,s}
\right),
$$
where $M_1=\mbox{const}>0$.
\end{lemma}

Proof follows from lemmas 2, 3 and 4.

\begin{lemma}
\label{lemma6}.
  For any function $f\in C^{2,s}(\bar {B}_{\rho _{j}})$
the following inequality holds:
$$
\left|
f
\right|_{(B_{\rho_{j}}) 1,s}
\leq
M_{2}
\left|
f
\right|_{(B_{\rho_{j}}) 2,s},
$$
where $M_2=\mbox{const}>0$.
\end{lemma}

{\bf Proof.}
We have the inequality (see [1]):
$$\sup_{u_{1},u_{2}\in B_{\rho_{j}}}
\frac{|\partial(f(u_{1}))-\partial(f(u_{2}))|}{|u_{1}-u_{2}|^{s}}\leq
$$
$$
\leq
M_{3}
\left(
\sum_{k=0}^{2}\sup_{B_{\rho_{j}}}|\partial^{k}(f)|
\right)^{s}
\left(
\sum_{k=0}^{1}\sup_{B_{\rho_{j}}}|\partial^{k}(f)|
\right)^{1-s}
\leq
M_{3}
\left(
\sum_{k=0}^{2}\sup_{B_{\rho_{j}}}|\partial^{k}(f)|
\right),
$$
where $M_3=\mbox{const}>0$.
Therefore we obtain the proof of lemma 6.

\section{Proof of theorem 1.}

By lemmas 5 and 6 we have the inequality:
$$
|f\circ h^{-1}_{i}(x^{1},...,x^{n-1},\omega(x^{1},...,x^{n-1}))|
_{(|(x^{1},...,x^{n-1})|\leq R) 2,s}
\leq
$$
$$
\leq
M_{4}|f\circ h^{-1}_{j}(u)|_{B_{\rho_{j}} 2,s},
$$
where $M_4=\mbox{const}>0$.
Since the ball $\bar {K}_{r_{i}}$ is compact, then
for some $j$ we obtain:
$$
\left|
f
\right|_{\partial K_{r_{i}}2,s}
\leq
M_{5}
\left|
f
\right|_{B_{\rho_{j}}2,s} \eqno(4),
$$
where $M_5=\mbox{const}>0$.

Since function $f$ is a solution of class $C^{2,s}(F)$
of the problem $A f =\gamma$ then we have the inequality (see [1]):
$$
|f|_{B_{\rho_{j}}2,s}\leq M(B_{\rho_{j}})
\left(
|\gamma|_{K_{r_{j}}0,s}+\max_{K_{r_{j}}}|f|
\right)
,
$$
where $\bar{B}_{\rho_{j}}\subset K_{r_{j}}$.

We have the inequality (see [1]):
$$\max_{K_{r_{j}}}|f|
\leq
\max
\left(
\max_{\partial K_{r_{j}}}|f|;
\max_{K_{r_{j}}}\Bigl|\frac{\gamma}{c}\Bigr|
\right).
$$

Since $F$ is compact then there exists point $x_{0}\in F$ such that
$|f(x_{0})|=\max_{F}|f|$.
Therefore there exists neighborhood $U_{l}$ such that $x_{0}\in U_{l}$.
Hence, we obtain:
$$
\max_{K_{r_{j}}}|f|\leq \max_{K_{r_{l}}}|f|\leq
\max_{K_{r_{l}}}\Bigl|\frac{\gamma}{c}\Bigr|.
\eqno(5)
$$

We have the estimate (see [1]):
$$
|f|_{K_{r_{i}}2,s}\leq M_{6}
\left(
|\gamma|_{K_{r_{i}}0,s}+\max_{K_{r_{i}}}|f|+|f|_{\partial K_{r_{i}}2,s}
\right),
$$
where $M_6=\mbox{const}>0.$
Using inequalities (4) and (5),
we finish the proof of theorem 1.

{\bf References}
\begin{enumerate}
\item
O.A. Ladyjenskaya, N.N. Uraltseva. Linear and quasilinear equations of
elliptic type. M: Nauka, 1973.
\end{enumerate}

\end{document}